# Lomax distribution and asymptotical ML estimations based on record values for probability density function and cumulative distribution function

**Saman Hosseini[1], Dler Hussein Kadir[2], Kostas Triantafyllopoulos[3]**

## Abstract

Here in this paper, it is tried to obtain and compare the ML estimations based on upper record values and a random sample. In continue, some theorems have been proven about the behavior of these estimations asymptotically.

## Introduction

Lomax distribution is one of the essential distributions in the applied statistics. It has a vast application in different fields of science, from actuarial sciences to econometrics and industry. Lomax distribution is, in fact, a shifted Pareto distribution. Mathematically it is shifted to change the starting point of the domain to zero. In this paper, variable X is distributed as Lomax distribution if the probability density function (pdf) is as follow:

$$f(x;\theta)=\frac{1}{\theta}(1+x)^{-(\frac{1}{\theta}+1)};x\geq 0,\theta\geq 0. \tag{1}$$

Consequently, the cumulative distribution function (CDF) is obtained as:

$$F(x;\theta)=1-(1+x)^{-\frac{1}{\theta}};x\geq 0,\theta\geq 0. \tag{2}$$

Nasiri and Hosseini (2012) have written a paper related to statistical inferences based on upper record values for Lomax distribution. Here, in this paper, it is tried to obtain maximum likelihood estimation (MLE) for the pdf and CDF of Lomax distribution based on a random sample and values

[1] Correspondent author, he is a lecturer at department of computer science, Cihan University-Erbil, Kurdistan Region, Iraq. Email them at Saman.hosseini@cihanuniversity.ac.iq and s.hosseini.stat@gmail.com
[2] Dler Hussein Kadir is a lecturer at department of statistics, Salahaddin University-Erbil, Kurdistan Region, Iraq
[3] Kostas Triantafyllopoulos is a professor at department of statistics, university od Sheffield, England

of the upper record. Additionally, it has been proven, these estimations based on record values in some cased are efficient and suitable via different theorems.

Presuming $\underset{i=1,\ldots,n}{X_i}$ is a random sample of size n from a given pdf like $f(x;\theta)$ which has m (m<n) values of upper records, the joint probability density function is obtained as below (Chandler 1952)

$$\underset{(R_1,\ldots,R_m)}{f}(r_1,\ldots r_m;\theta) = f(r_m;\theta)\{\prod_{i=1}^{m-1} h(r_i;\theta).$$

Where $h(r_i;\theta) = \dfrac{f(r_i;\theta)}{1-F(r_i;\theta)}$.

Considering (1) and (2) joint probability density function for m values of upper records from Lomax distribution is easily obtained as follow:

$$\underset{(R_1,\ldots,R_m)}{f}(r_1,\ldots r_m;\theta) = \theta^m \frac{(1+r_m)^{-\frac{1}{\theta}}}{\prod_{i=1}^{m}(1+r_i)} \tag{3}$$

# Maximum Likelihood Estimations (MLE) for the unknown parameter $\theta$ based on a random sample and record values.

Using (3), MLE based on m values of the upper record is obtained as (in this paper is shown by $\hat{\theta}_{\mathrm{Re}cords;m}$): $\hat{\theta}_{\mathrm{Re}cords;m} = \dfrac{\ln(1+R_m)}{m}$. (Arnold, Balakrishnan and Nagaraja, 1998)

On the other hand, if $\underset{i=1,\ldots,n}{X_i}$ is a random sample of size n from Lomax distribution, the ML estimation based on this sample for the unknown parameter $\theta$ is obtained easily:

$$\underset{(X_1,\ldots,X_n)}{f}(x_1,\ldots,x_n;\theta) = \prod_{i=1}^{n} f_{X_i}(x_i;\theta) = \prod_{i=1}^{n}\frac{1}{\theta}(1+x_i)^{-(\frac{1}{\theta}+1)} = \theta^{-n}\{\prod_{i=1}^{n}(1+x_i)\}^{-(\frac{1}{\theta}+1)}.$$

Then

$$L = \ln(\underset{(X_1,\ldots,X_n)}{f}(x_1,\ldots,x_n;\theta)) = -n\ln(\theta) - (\frac{1}{\theta}+1)\sum_{i=1}^{n}\ln(1+x_i).$$

By differentiation and some simple arithmetical operations, the ML estimation based on a random sample of size n ($\theta_{Sample,n}$) is obtained.

$$\frac{\partial L}{\partial \theta}=0 \Rightarrow \frac{-n}{\theta}+\frac{\sum_{i=1}^{n}\ln(1+x_i).}{\theta^2}=0 \Rightarrow \theta=\frac{\sum_{i=1}^{n}\ln(1+x_i).}{n} \Rightarrow \hat{\theta}_{sample,n}=\frac{\sum_{i=1}^{n}\ln(1+X_i)}{n}. \qquad (4)$$

Above value is an ML estimation, because the sign of the second derivation for this point is negative:

$$\frac{\partial L}{\partial \theta}=\frac{-n}{\theta}+\frac{\sum_{i=1}^{n}\ln(1+x_i)}{\theta^2} \Rightarrow \frac{\partial^2 L}{\partial \theta^2}=\frac{n}{\theta^2}-2\frac{\sum_{i=1}^{n}\ln(1+x_i)}{\theta^3}$$

$$\Rightarrow \left.\frac{\partial^2 L}{\partial \theta^2}\right|_{\hat{\theta}_{sample,n}=\frac{\sum_{i=1}^{n}\ln(1+X_i)}{n}}=\frac{n}{\left(\frac{\sum_{i=1}^{n}\ln(1+X_i)}{n}\right)^2}-2\frac{\sum_{i=1}^{n}\ln(1+x_i)}{\left(\frac{\sum_{i=1}^{n}\ln(1+X_i)}{n}\right)^3}$$

$$\Rightarrow \left.\frac{\partial^2 L}{\partial \theta^2}\right|_{\hat{\theta}_{sample,n}=\frac{\sum_{i=1}^{n}\ln(1+X_i)}{n}}=\frac{n^3}{(\sum_{i=1}^{n}\ln(1+X_i))^2}-2\frac{n^3}{(\sum_{i=1}^{n}\ln(1+X_i))^2}=-\frac{n^3}{(\sum_{i=1}^{n}\ln(1+X_i))^2}<0.$$

The following theorem explains one of the most primary reasons for using ML estimations based on upper record values for the unknown parameter of Lomax. In other words, the following theorem gives some prestige to use MLE based on upper records for this distribution.

## Theorem 1.

If $\underset{i=1,\dots,n}{X_i}$ is a random sample of size n from Lomax distribution and $(R_1,\dots,R_n)$ are n values of the upper record from Lomax distribution and $g(t)$ is an on to one real function, then:

I- $g(\hat{\theta}_{\text{Re}cords;n})$ is equally distributed with $g(\hat{\theta}_{Sample;n})$. Mathematically it means

$$g(\hat{\theta}_{\text{Re}cords;n}) \overset{distribution}{=} g(\hat{\theta}_{Sample;n})$$

II- $MSE[g(\hat{\theta}_{\text{Re}cords;n})]=MSE[g(\hat{\theta}_{Sample;n})]$.

## Proof.

For the first part of the theorem the distribution of $g(\hat{\theta}_{\text{Re}cords;n})$ is obtained:

The probability density function of $R_n$ is obtained as: $f_{R_n}(r_n)=f(r)\frac{[-\ln(1-F(r))]^{n-1}}{(n-1)!}$.

**(**Arnold, Balakrishnan and Nagaraja, 1998**)**

On the other hand, $f(x;\theta)=\frac{1}{\theta}(1+x)^{-(\frac{1}{\theta}+1)}$ and $F(x;\theta)=1-(1+x)^{-\frac{1}{\theta}}$ then

$$f_{R_n}(r_n)=\frac{1}{\theta}(1+r_n)^{-(\frac{1}{\theta}+1)}\frac{[-\ln(1-\{1-(1+r_n)^{-\frac{1}{\theta}}\})]^{n-1}}{(n-1)!}=\frac{1}{\theta^n}(1+r_n)^{-(\frac{1}{\theta}+1)}\frac{[\ln(1+r_n)]^{n-1}}{(n-1)!} \tag{5}$$

If $T=\ln(1+R_n)$ and considering (5)

$$F_T(t)=P(T\le t)=P(\ln(1+R_n)\le t)=P(R_n\le e^t-1)=F_{R_n}(e^t-1)$$

consequently

$$F_T(t)=F_{R_n}(e^t-1)\Rightarrow\frac{dF_T(t)}{dt}=\frac{F_{R_n}(e^t-1)}{dt}\Rightarrow f_T(t)=e^t f_{R_n}(e^t-1)$$

by substitution $e^t-1$ in (5) and multiplying obtained function by $e^t$ we have:

$$f_T(t)=e^t\frac{1}{\theta^n}(1+e^t-1)^{-(\frac{1}{\theta}+1)}\frac{[\ln(1+e^t-1)]^{n-1}}{(n-1)!}=\frac{t^{n-1}\exp\{-\frac{t}{\theta}\}}{\theta^n\Gamma(n)}.$$

It means $T=\ln(1+R_n)$ is distributed as $gamma(n,\theta)$, and consequently

$$\ln(1+R_n)\overset{distribution}{=}T\sim Gamma(n,\theta)$$

Consequently

$$\hat{\theta}_{\mathrm{Re}cords;n}=\frac{\ln(1+R_n)}{n}\overset{distribution}{=}\frac{T}{n}\Rightarrow g(\hat{\theta}_{\mathrm{Re}cords;n})\overset{distribution}{=}g(\frac{T}{n}) \tag{6}$$

Now the distribution of $g(\hat{\theta}_{Sample;n})$ is obtained. As it was mentioned and obtained in (4)

$\hat{\theta}_{sample,n}=\frac{\sum_{i=1}^{n}\ln(1+X_i)}{n}$, on the other hand, $\underset{i=1,\dots,n}{X_i}$ are independent and distributed identically from Lomax distribution. Then the following can be obtained easily:

$$P(\ln(1+X_i)\le x)=P(1+X_i\le e^x)=P(X_i\le e^x-1))\Rightarrow F_{\ln(1+X_i)}(x)=F_{X_i}(e^x-1).$$

Consequently:

$$f_{\ln(1+X_i)}(x)=e^x f_{X_i}(e^x-1)$$

Substituting $(e^x - 1)$ in (1) and multiplying it by $e^x$ the following is concluded:

$$f_{\ln(1+X_i)}(x) = \frac{e^x}{\theta}(1+(e^x - 1))^{-(\frac{1}{\theta}+1)} = \frac{1}{\theta}e^{-\frac{x}{\theta}}$$

It means $\ln(1+X_i)$ is exponentially distributed, and its parameter is $\theta$ ( $\ln(1+X_i) \overset{distribution}{=} \exp(\theta)$ ). Additionally, as was mentioned before $\underset{i=1,\dots,n}{X_i}$ are independent, and that's why each one-to-one function of them is independent as well. Then $\underset{i=1,\dots,n}{\ln(1+X_i)}$ are independent; on the other hand, each $\ln(1+X_i)$ is exponentially distributed. Considering the above explanation, it is easy to obtain the distribution of $\frac{\sum_{i=1}^{n}\ln(1+X_i)}{n}$.

$$\ln(1+X_i) \overset{distribution}{=} \exp(\theta) \Rightarrow \sum_{i=1}^{n}\ln(1+X_i) \overset{distribution}{=} T \sim gamma(n,\theta) \Rightarrow \frac{\sum_{i=1}^{n}\ln(1+X_i)}{n} \overset{distribution}{=} \frac{T}{n}$$

It means

$$\hat{\theta}_{Sample;n} \overset{distribution}{=} \frac{T}{n} \Rightarrow g(\hat{\theta}_{Sample;n}) \overset{distribution}{=} g(\frac{T}{n}) \quad (7)$$

By (6) and (7), the first part of the theorem is proven. The second part of the theorem is a direct conclusion of the first part.

## Asymptotically Properties of ML Estimations based on records for PDF and CDF of Lomax distribution.

Considering the invariance property of ML estimations, ML estimations of PDF and CDF are easily obtained as follow:

$$\hat{f}_{(Records;m)}(x) = \frac{m}{\ln(1+R_m)}(1+x)^{-(\frac{m}{\ln(1+R_m)}+1)}$$

and

$$\hat{F}_{(Records;m)}(x) = 1-(1+x)^{-(\frac{m}{\ln(1+R_m)})}$$

The following theorems show that these estimations are useful if the number of record values exceeds a limit.

## Theorem 2.

If $\underset{i=1,\dots,n}{X_i}$ is a random sample from Lomax distribution and contains m values of upper record $(R_1,\dots,R_m)$ then:

$$I-E[\hat{f}_{(\mathrm{Re}cords;m)}(x)]=(\frac{m}{\theta(1+x)})\sum_{i=0}^{m-1}\frac{\Gamma(m-i-1)}{\Gamma(m)\Gamma(i+1)}(\frac{-m\ln(1+x)}{\theta})^i$$

$$II-E[\hat{F}_{(\mathrm{Re}cords;m)}(x)]=1-\sum_{i=0}^{m}\frac{\Gamma(m-i)}{\Gamma(m)\Gamma(i+1)}(\frac{-m\ln(1+x)}{\theta})^i$$

**Proof.** It is known by (6) that $\hat{f}_{(\mathrm{Re}cords;m)}(x)=\frac{m}{\ln(1+R_m)}(1+x)^{-(\frac{m}{\ln(1+R_m)}+1)}\overset{distribution}{=}(\frac{m}{T})(1+x)^{-(\frac{m}{T}+1)}$

where $T\sim gamma(m,\theta)$ then

$$E[\hat{f}_{(\mathrm{Re}cords;m)}(x)]=\int_0^{+\infty}(\frac{m}{t})(1+x)^{-(\frac{m}{t}+1)}\frac{t^{m-1}}{\Gamma(m)\theta^m}e^{-\frac{t}{\theta}}dt=\frac{m}{\Gamma(m)\theta^m}\int_0^{+\infty}(1+x)^{-(\frac{m}{t}+1)}t^{m-2}e^{-\frac{t}{\theta}}dt$$

Considering $(1+x)^{-\frac{m}{t}}=\exp(-\frac{m}{t}\ln(1+x))$ the following is concluded:

$$E[\hat{f}_{(\mathrm{Re}cords;m)}(x)]=\frac{m}{\Gamma(m)\theta^m}\int_0^{+\infty}(1+x)^{-(\frac{m}{t}+1)}t^{m-2}e^{-\frac{t}{\theta}}dt=\frac{m}{\Gamma(m)\theta^m(1+x)}\int_0^{+\infty}(1+x)^{-\frac{m}{t}}t^{m-2}e^{-\frac{t}{\theta}}dt$$

$$=\frac{m}{\Gamma(m)\theta^m(1+x)}\int_0^{+\infty}e^{-\frac{m}{t}\ln(1+x)}t^{m-2}e^{-\frac{t}{\theta}}dt=\frac{m}{\Gamma(m)\theta^m(1+x)}\int_0^{+\infty}\sum_{i=0}^{+\infty}\frac{(-m\ln(1+x))^i}{t^i i!}t^{m-2}e^{-\frac{t}{\theta}}dt$$

$$=\frac{m}{\Gamma(m)\theta^m(1+x)}\int_0^{+\infty}\sum_{i=0}^{+\infty}\frac{(-m\ln(1+x))^i}{\Gamma(i+1)}t^{m-i-2}e^{-\frac{t}{\theta}}dt$$

$$=\frac{m}{\Gamma(m)\theta^m(1+x)}\int_0^{+\infty}\sum_{i=0}^{+\infty}\frac{(-m\ln(1+x))^i}{\Gamma(i+1)}t^{m-i-2}e^{-\frac{t}{\theta}}dt=\frac{m}{\Gamma(m)\theta^m(1+x)}\sum_{i=0}^{+\infty}\frac{(-m\ln(1+x))^i}{\Gamma(i+1)}\int_0^{+\infty}t^{m-i-2}e^{-\frac{t}{\theta}}dt$$

By the simple transformation $\frac{t}{\theta}=z$ and considering $t=\theta z\Rightarrow dt=\theta dz$ the following relations are obtained

$$E[\hat{f}_{(\text{Re}cords;m)}(x)] = \frac{m}{\Gamma(m)\theta^m(1+x)}\sum_{i=0}^{+\infty}\frac{(-m\ln(1+x))^i}{\Gamma(i+1)}\int_0^{+\infty} t^{m-i-2}e^{-\frac{t}{\theta}}dt =$$

$$\frac{m}{\Gamma(m)\theta^m(1+x)}\sum_{i=0}^{+\infty}\frac{(-m\ln(1+x))^i}{\Gamma(i+1)}\Gamma(m-i-1)\theta^{m-i-1}$$

$$= \frac{m}{\theta(1+x)}\sum_{i=0}^{+\infty}\frac{\Gamma(m-i-1)}{\Gamma(i+1)\Gamma(m)}\left(\frac{-m\ln(1+x)}{\theta}\right)^i \quad (8)$$

Finally, the first part is proven. The second part has similar proof.

## Theorem 3.

if $\underset{i=1,\ldots,n}{X_i}$ is a random sample from Lomax distribution and contains m values of upper record $(R_1,\ldots,R_m)$ , $\hat{f}_{(\text{Re}cords;m)}(x)$ and $\hat{F}_{(\text{Re}cords;m)}(x)$ are asymptotically unbiased for PDF and CDF of Lomax distribution. It means:

$$I - \lim_{m\to+\infty} E[\hat{f}_{(\text{Re}cords;m)}(x)] = f(x;\theta) = \frac{1}{\theta}(1+x)^{-(\frac{1}{\theta}+1)}$$

$$II - \lim_{m\to+\infty} E[\hat{F}_{(\text{Re}cords;m)}(x)] = F(x;\theta) = F(x;\theta) = 1-(1+x)^{-\frac{1}{\theta}}$$

## Proof.

First, we need to state the following lemma:

## Lemma1.

$$\underset{n\to\infty}{Lim}\ \frac{\Gamma(n-i-1)n^{i+1}}{\Gamma(n)} = 1.$$

Gore, Hosseini, & Nasiri (2017) proved this lemma.

The first part is proven; the second part has similar proof. Considering (7)

$$\lim_{m\to+\infty} E[\hat{f}_{(\text{Re}cords;m)}(x)] = \lim_{m\to+\infty}\frac{m}{\theta(1+x)}\sum_{i=0}^{+\infty}\frac{\Gamma(m-i-1)}{\Gamma(i+1)\Gamma(m)}\left(\frac{-m\ln(1+x)}{\theta}\right)^i$$

Then it means:

$$\lim_{m\to+\infty} E[\hat{f}_{(\mathrm{Re}cords;m)}(x)] = \lim_{m\to+\infty} \frac{}{\theta(1+x)} \sum_{i=0}^{+\infty} \frac{\Gamma(m-i-1)}{\Gamma(i+1)\Gamma(m)} (\frac{-m\ln(1+x)}{\theta})^i$$

$$= \lim_{m\to+\infty} \frac{1}{\theta(1+x)} \sum_{i=0}^{+\infty} \frac{\Gamma(m-i-1)m^{i+1}}{\Gamma(m)} \frac{1}{\Gamma(i+1)} (\frac{-\ln(1+x)}{\theta})^i$$

$$= \frac{1}{\theta(1+x)} \sum_{i=0}^{+\infty} (\underbrace{\lim_{m\to+\infty} \frac{\Gamma(m-i-1)m^{i+1}}{\Gamma(m)}}_{lemma1.}) \frac{1}{\Gamma(i+1)} (\frac{-\ln(1+x)}{\theta})^i$$

$$= \frac{1}{\theta(1+x)} \sum_{i=0}^{+\infty} (1) \frac{1}{\Gamma(i+1)} (\frac{-\ln(1+x)}{\theta})^i = \frac{1}{\theta(1+x)} \sum_{i=0}^{+\infty} \frac{(\frac{-\ln(1+x)}{\theta})^i}{i!} = \frac{1}{\theta(1+x)} \exp(\ln(1+x)^{-\frac{1}{\theta}}) = \frac{1}{\theta}(1+x)^{-\frac{1}{\theta}-1}$$

Finally, $\lim_{m\to+\infty} E[\hat{f}_{(\mathrm{Re}cords;m)}(x)] = f(x;\theta) = \frac{1}{\theta}(1+x)^{-\frac{1}{\theta}-1}$.

After theorem 3 the following theorem is interesting; it gives us a simple and attractive relation between estimations of PDF ( $E[\hat{f}_{(\mathrm{Re}cords;m)}(x)]$ ) and CDF ( $E[\hat{F}_{(\mathrm{Re}cords;m)}(x)]$ ) based on records.

## Theorem 4.

If a given sample like $\underset{i=1,\dots,n}{X_i}$ from Lomax distribution contains m values of upper record $(R_1,\dots,R_m)$ , the following relation is always true asymptotically:

$I) 1 - E[\hat{F}_{(\mathrm{Re}cords;m)}(x)] = \theta(1+x) E[\hat{f}_{(\mathrm{Re}cords;m)}(x)]$

$II) E[\hat{f}_{(\mathrm{Re}cords;m)}(x)] < \frac{1}{\theta}$

## Proof.

Considering relation (1), (2), and theorem 3, the following is easily obtained:

$$\left.\begin{aligned} f(x;\theta) &= \frac{1}{\theta}(1+x)^{-(\frac{1}{\theta}+1)} \\ F(x;\theta) &= 1-(1+x)^{-\frac{1}{\theta}} \end{aligned}\right\} \Rightarrow \frac{1-F(x;\theta)}{f(x;\theta)} = \frac{(1+x)^{-\frac{1}{\theta}}}{\frac{1}{\theta}(1+x)^{-(\frac{1}{\theta}+1)}} = \theta(1+x) \qquad (9)$$

By theorem 3 and substituting $\lim_{m\to+\infty} E[\hat{F}_{(\mathrm{Re}cords;m)}(x)] = F(x;\theta)$ and $\lim_{m\to+\infty} E[\hat{f}_{(\mathrm{Re}cords;m)}(x)] = f(x;\theta)$ in (9) we have

$$\frac{1 - \lim_{m\to+\infty} E[\hat{F}_{(\mathrm{Re}cords;m)}(x)]}{\lim_{m\to+\infty} E[\hat{f}_{(\mathrm{Re}cords;m)}(x)]} = \theta(1+x)$$

it means $\exists N_0 \in \mathbb{N}, n \geq N_0$ such that $\dfrac{1-E[\hat{F}_{(\mathrm{Re}cords;m)}(x)]}{E[\hat{f}_{(\mathrm{Re}cords;m)}(x)]} = \theta(1+x)$.

Considering relation (1), it is easy to find the maximum value of PDF from Lomax:

$$f(x;\theta) = \frac{1}{\theta}(1+x)^{-(\frac{1}{\theta}+1)} \Rightarrow \frac{df(x;\theta)}{dx} = -\frac{1}{\theta}(\frac{1}{\theta}+1)(1+x)^{-(\frac{1}{\theta}+2)}$$

The derivative is negative; that's why it is a decreasing function of x. Therefore the maximum value is obtained from the beginning point of the domain ( $x=0$ ), and the maximum value of PDF is

$$\max f(x;\theta) = f(0;\theta) = \frac{1}{\theta}(1+0)^{-(\frac{1}{\theta}+1)} = \frac{1}{\theta}$$

It means

$$\forall x \geq 0;\ f(x;\theta) \leq \frac{1}{\theta} \text{ and } \forall x \geq 0;\ \lim_{m\to+\infty} E[\hat{f}_{(\mathrm{Re}cords;m)}(x)] \leq \frac{1}{\theta}$$

Consequently $\exists N_0 \in \mathbb{N}, n \geq N_0$

$$E[\hat{f}_{(\mathrm{Re}cords;m)}(x)] < \frac{1}{\theta}$$

## Theorem 5.

$$I - MSE(\hat{f}_{(\mathrm{Re}cords,m)}) = (\frac{m}{\theta(1+x)})^2 \sum_{i=0}^{m-2} \frac{\Gamma(m-i-2)}{\Gamma(m)\Gamma(i+1)}(-2\frac{m\ln(1+x)}{\theta})^i -$$
$$\frac{2m}{\theta^2}(1+x)^{-(\frac{1}{\theta}+2)} \sum_{i=0}^{m} \frac{\Gamma(m-i)}{\Gamma(m)\Gamma(i+1)}(-\frac{m\ln(1+x)}{\theta})^i + \frac{1}{\theta^2}(1+x)^{-2(\frac{1}{\theta}+1)}$$

$$II - MSE(\hat{F}_{(\mathrm{Re}cords,m)}) = \sum_{i=0}^{m} \frac{\Gamma(m-i)}{\Gamma(m)\Gamma(i+1)}(-2\frac{m\ln(1+x)}{\theta})^i -$$
$$2(1+x)^{-\frac{1}{\theta}} \sum_{i=0}^{m} \frac{\Gamma(m-i)}{\Gamma(m)\Gamma(i+1)}(-\frac{m\ln(1+x)}{\theta})^i + (1+x)^{-\frac{2}{\theta}}.$$

## Proof.

It is clear that for a variable like T used to estimate $A(\theta)$, the Mean Square Error (MSE) is defined as: $MSE(T) = E[T - A(\theta)]^2 = E[T^2] - 2A(\theta)E[T] + (A(\theta))^2$. Based on this definition, MSE is calculated for $\hat{f}_{(\mathrm{Re}cords,m)}$ as follow:

$$MSE(\hat{f}_{(\mathrm{Re}cords,m)}) = E[(\hat{f}_{(\mathrm{Re}cords,m)})^2] - 2f(x)E[\hat{f}_{(\mathrm{Re}cords,m)}] + (f(x))^2 \qquad (10)$$

Considering (6) $E[(\hat{f}_{(\mathrm{Re}cords,m)})^2]=E[\frac{m^2}{T^2}(1+x)^{-2(\frac{m}{T}+1)}]$ in which $T$ is distributed as $gamma(m,\theta)$, that's why:

$$E[(\hat{f}_{(\mathrm{Re}cords,m)})^2]=E[\frac{m^2}{T^2}(1+x)^{-2(\frac{m}{T}+1)}]=\int_0^{+\infty}\frac{m^2}{t^2}(1+x)^{-2(\frac{m}{t}+1)}\frac{t^{m-1}}{\Gamma(m)\theta^m}e^{-\frac{t}{\theta}}dt=$$

$$\frac{m^2}{(1+x)^2\Gamma(m)\theta^m}\int_0^{+\infty}(1+x)^{-2\frac{m}{t}}t^{m-3}e^{-\frac{t}{\theta}}dt$$

For ease, it is better to consider $\gamma(m,x,\theta)=\frac{m^2}{(1+x)^2\Gamma(m)\theta^m}$.

Now

$$E[(\hat{f}_{(\mathrm{Re}cords,m)})^2]=\lambda(m,x,\theta)\int_0^{+\infty}(1+x)^{-2\frac{m}{t}}t^{m-3}e^{-\frac{t}{\theta}}dt=\lambda(m,x,\theta)\int_0^{+\infty}\mathrm{e}^{(-2\frac{m}{t}\ln(1+x))}t^{m-3}e^{-\frac{t}{\theta}}dt$$

$$=\lambda(m,x,\theta)\int_0^{+\infty}\sum_{i=0}^{+\infty}\frac{\{-2m\ln(1+x)\}^i}{i!t^i}t^{m-3}e^{-\frac{t}{\theta}}dt=\lambda(m,x,\theta)\int_0^{+\infty}\sum_{i=0}^{+\infty}\frac{\{-2m\ln(1+x)\}^i}{\Gamma(i+1)}t^{m-i-3}e^{-\frac{t}{\theta}}dt$$

If $A(x)=\ln(1+x)$, then:

$$E[(\hat{f}_{(\mathrm{Re}cords,m)})^2]=\lambda(m,x,\theta)\int_0^{+\infty}\sum_{i=0}^{+\infty}\frac{\{-2mA(x)\}^i}{\Gamma(i+1)}t^{m-i-3}e^{-\frac{t}{\theta}}dt=$$

$$\lambda(m,x,\theta)\sum_{i=0}^{+\infty}\frac{\{-2mA(x)\}^i}{\Gamma(i+1)}\int_0^{+\infty}t^{m-i-3}e^{-\frac{t}{\theta}}dt$$

Substituting $\frac{t}{\theta}=z$ in the last integration

$$E[(\hat{f}_{(\mathrm{Re}cords,m)})^2]=\lambda(m,x,\theta)\sum_{i=0}^{+\infty}\frac{\{-2mA(x)\}^i}{\Gamma(i+1)}\int_0^{+\infty}\theta^{m-i-2}z^{m-i-3}e^{-z}dz=$$

$$\lambda(m,x,\theta)\sum_{i=0}^{+\infty}\frac{\{-2mA(x)\}^i}{\Gamma(i+1)}\theta^{m-i-2}\Gamma(m-i-2)$$

Finally;

$$E[(\hat{f}_{(\mathrm{Re}cords,m)})^2]=(\frac{m}{\theta(1+x)})^2\sum_{i=0}^{+\infty}\frac{\Gamma(m-i-2)}{\Gamma(i+1)\Gamma(m)}\{-2\frac{m\ln(1+x)}{\theta}\}^i \qquad (11)$$

Considering (8), (10), and (11), the following is concluded:

$$MSE(\hat{f}_{(\mathrm{Re}cords,m)}) = (\frac{m}{\theta(1+x)})^2 \sum_{i=0}^{m-2} \frac{\Gamma(m-i-2)}{\Gamma(m)\Gamma(i+1)} (-2\frac{m\ln(1+x)}{\theta})^i -$$

$$\frac{2m}{\theta^2}(1+x)^{-(\frac{1}{\theta}+2)} \sum_{i=0}^{m} \frac{\Gamma(m-i)}{\Gamma(m)\Gamma(i+1)} (-\frac{m\ln(1+x)}{\theta})^i + \frac{1}{\theta^2}(1+x)^{-2(\frac{1}{\theta}+1)}$$

And the first part of the theorem is proven. The proof for the second part is similar.

## Theorem 6.

$I - \hat{f}_{(\mathrm{Re}cords,m)} \xrightarrow{P} f(x;\theta)$

$II - \hat{F}_{(\mathrm{Re}cords,m)} \xrightarrow{P} F(x;\theta)$

$III - \hat{\theta}_{(\mathrm{Re}cords,m)} \xrightarrow{P} \theta$

## Proof.

Based on Markov theorem, the following statement is clear:

$$P(\left| \hat{f}_{(\mathrm{Re}cords,m)} - E(\hat{f}_{(\mathrm{Re}cords,m)}) \right| \ge \in) \le {Var(\hat{f}_{(\mathrm{Re}cords,m)})} \big/ {\in}$$

On the other hand, by (8) and (11) we have

$$E[(\hat{f}_{(\mathrm{Re}cords,m)})^2] = (\frac{m}{\theta(1+x)})^2 \sum_{i=0}^{+\infty} \frac{\Gamma(m-i-2)}{\Gamma(i+1)\Gamma(m)} \{-2\frac{m\ln(1+x)}{\theta}\}^i \text{ and}$$

$$E[\hat{f}_{(\mathrm{Re}cords,m)}] = \frac{m}{\theta(1+x)} \sum_{i=0}^{+\infty} \frac{\Gamma(m-i-1)}{\Gamma(i+1)\Gamma(m)} (\frac{-m\ln(1+x)}{\theta})^i \text{, therefore:}$$

$Var(\hat{f}_{(\mathrm{Re}cords,m)}) = E[(\hat{f}_{(\mathrm{Re}cords,m)})^2] - E^2[\hat{f}_{(\mathrm{Re}cords,m)}]$ and consequently:

$$Var(\hat{f}_{(\mathrm{Re}cords,m)}) = (\frac{m}{\theta(1+x)})^2 \sum_{i=0}^{m-2} \frac{\Gamma(m-i-2)}{\Gamma(i+1)\Gamma(m)} \{-2\frac{m\ln(1+x)}{\theta}\}^i -$$

$$(\frac{m}{\theta(1+x)} \sum_{i=0}^{m-1} \frac{\Gamma(m-i-1)}{\Gamma(i+1)\Gamma(m)} (\frac{-m\ln(1+x)}{\theta})^i)^2$$

By limitation:

$$\lim_{m\to+\infty} Var(\hat{f}_{(\mathrm{Re}cords,m)}) = \lim_{m\to+\infty}\{(\frac{m}{\theta(1+x)})^2\sum_{i=0}^{m-2}\frac{\Gamma(m-i-2)}{\Gamma(i+1)\Gamma(m)}\{-2\frac{m\ln(1+x)}{\theta}\}^i - (E[\hat{f}_{(\mathrm{Re}cords,m)}])^2\}$$

$$\lim_{m\to+\infty} Var(\hat{f}_{(\mathrm{Re}cords,m)}) = \lim_{m\to+\infty}\{(\frac{m}{\theta(1+x)})^2\sum_{i=0}^{m-2}\frac{\Gamma(m-i-2)}{\Gamma(i+1)\Gamma(m)}\{-2\frac{m\ln(1+x)}{\theta}\}^i\} - \lim_{m\to+\infty}\{(E[\hat{f}_{(\mathrm{Re}cords,m)}])^2\}$$

And during the first part of Theorem 3, it was discussed that $\lim_{m\to+\infty} E[\hat{f}_{(\mathrm{Re}cords,m)}] = f(x;\theta)$, then

$$\lim_{m\to+\infty} Var(\hat{f}_{(\mathrm{Re}cords,m)}) = \lim_{m\to+\infty}\{(\frac{m}{\theta(1+x)})^2\sum_{i=0}^{m-2}\frac{\Gamma(m-i-2)}{\Gamma(i+1)\Gamma(m)}\{-2\frac{m\ln(1+x)}{\theta}\}^i\} - \{f(x;\theta)\}^2$$

$$\lim_{m\to+\infty} Var(\hat{f}_{(\mathrm{Re}cords,m)}) = \lim_{m\to+\infty}\{(\frac{m}{\theta(1+x)})^2\sum_{i=0}^{m-2}\frac{\Gamma(m-i-2)m^{i+2}}{\Gamma(m)}\frac{1}{m^2\Gamma(i+1)}\{-2\frac{\ln(1+x)}{\theta}\}^i\} - \{f(x;\theta)\}^2$$

Considering lemma 1. The following is concluded:

$$\lim_{m\to+\infty} Var(\hat{f}_{(\mathrm{Re}cords,m)}) =$$

$$\lim_{m\to+\infty}\{(\frac{m}{\theta(1+x)})^2\frac{1}{m^2}\}\sum_{i=0}^{+\infty}\lim_{m\to+\infty}\{\frac{\Gamma(m-i-2)m^{i+2}}{\Gamma(m)}\}\frac{1}{\Gamma(i+1)}\{-2\frac{\ln(1+x)}{\theta}\}^i - \{f(x;\theta)\}^2$$

$$= \lim_{m\to+\infty}\{(\frac{1}{\theta(1+x)})^2\}\sum_{i=0}^{+\infty}\frac{1}{\Gamma(i+1)}\{-2\frac{\ln(1+x)}{\theta}\}^i - \{f(x;\theta)\}^2$$

$$= \{(\frac{1}{\theta(1+x)})^2\}\sum_{i=0}^{+\infty}\frac{1}{\Gamma(i+1)}\{-2\frac{\ln(1+x)}{\theta}\}^i - \{f(x;\theta)\}^2 = \frac{1}{\theta^2}(1+x)^{-2(\frac{1}{\theta}+1)} - \{f(x;\theta)\}^2 = 0$$

It means

$$\lim_{m\to+\infty} P(\left|\hat{f}_{(\mathrm{Re}cords,m)} - E(\hat{f}_{(\mathrm{Re}cords,m)})\right| \geq \in) \leq 0 \Rightarrow \lim_{m\to+\infty} P(\left|\hat{f}_{(\mathrm{Re}cords,m)} - E(\hat{f}_{(\mathrm{Re}cords,m)})\right| \geq \in) = 0$$

Consequently:

$\hat{f}_{(\mathrm{Re}cords,m)} - E(\hat{f}_{(\mathrm{Re}cords,m)}) \xrightarrow{P} 0$, or $\hat{f}_{(\mathrm{Re}cords,m)} \xrightarrow{P} E(\hat{f}_{(\mathrm{Re}cords,m)})$.

On the other hand, $\lim_{m\to+\infty} E[\hat{f}_{(\mathrm{Re}cords,m)}] = f(x;\theta)$ then:

$$\left.\begin{array}{l}\lim_{m\to+\infty} E[\hat{f}_{(\mathrm{Re}cords,m)}] = f(x;\theta) \\ \hat{f}_{(\mathrm{Re}cords,m)} - E(\hat{f}_{(\mathrm{Re}cords,m)}) \xrightarrow{P} 0\end{array}\right\}$$

$$\Rightarrow \hat{f}_{(\mathrm{Re}cords,m)} - E(\hat{f}_{(\mathrm{Re}cords,m)}) + \lim_{m\to+\infty} E[\hat{f}_{(\mathrm{Re}cords,m)}] \xrightarrow{P} 0 + f(x;\theta)$$

It means the first part of the theorem is proven. The proof for the second part is similar.

For the third part of the theorem, the following lemma is necessary to state.

**Lemma 2.**

If $Z_m \xrightarrow{P} z$ and $D:\mathbb{R}\rightarrow\mathbb{R}$ is a continuous function then $D(Z_m)\xrightarrow{P} D(z)$ (Billingsley, 1995)

Considering lemma 2 and the second part of the theorem 6:

If $D(t) = -\frac{\ln(1+x)}{\ln(1-t)}$ , on the other hand, based on the second part of theorem 6

$\hat{F}_{(Records,m)} \xrightarrow{P} F(x;\theta)$ then:

$$DoF\hat{F}_{(Records,m)} \xrightarrow{P} DoF(x;\theta)$$

Consequently:

$$\left.\begin{array}{r} DoF\hat{F}_{(Records,m)} = -\dfrac{\ln(1+x)}{\ln(1-(1-(1+x)^{\hat{\theta}_{(Records,m)}}))} = \hat{\theta}_{(Records,m)} \\ DoF(x;\theta) = -\dfrac{\ln(1+x)}{\ln(1-(1-(1+x)^{-\frac{1}{\theta}}))} = \theta \end{array}\right\} \Rightarrow \hat{\theta}_{(Records,m)} \xrightarrow{P} \theta .$$

## Conclusion.

In this paper, ML estimations for the unknown parameter, PDF, and CDF of Lomax distribution were obtained. Afterward, some exciting relations between estimations based on a random sample of size n and n values of record value were obtained and proven during some theorems. Next by some theorems, it was shown that at least in asymptotical situations, it worth to use estimation based on records values. It means based on proven theorems in this paper; asymptotical estimations are unbiased and convergent in probability to the related parameters, which generally means in the case of having many values of record, estimation based on them is enough trustable compared with estimations based on a random sample.

## References.

1- Nasiri, P., & Hosseini, S. (2012). Statistical inferences for Lomax distribution based on record values (Bayesian and classical). Journal of Modern Applied Statistical Methods, 11(1), 15.
2- Chandler, K. N. (1952). The distribution and frequency of record values. Journal of the Royal Statistical Society. Series B (Methodological), 220-228.
3- B. C. Arnold, N. Balakrishnan, H.N. Nagaraja,"records,"John Wiley and Sons, Canada1998.
4- arXiv: 1710.10690 [math. ST]
5- Billingsley, P. (1995). Probability and Measure, ser. Probability and Mathematical Statistics. New York: Wiley.